\documentclass[english]{iopart}

\usepackage{graphicx}
\usepackage{url}
\usepackage{iopams} 

\expandafter\let\csname equation*\endcsname\relax
\expandafter\let\csname endequation*\endcsname\relax

\usepackage{amssymb}
\usepackage[latin1]{inputenc}
\usepackage{amsmath}
\usepackage{epsfig}
\newcounter{fig}
\begin{document}

    \title[Christol's conjecture]
{\Large On Christol's conjecture}

\vskip .3cm 

\author{Y. Abdelaziz$^\dag$, C. Koutschan$^\P$,
J-M. Maillard$^\dag$}

\address{$^\dag$ LPTMC, UMR 7600 CNRS, 
Universit\'e Pierre et Marie Curie, Sorbonne Universit\'e, 
Tour 23, 5\`eme \'etage, case 121, 
 4 Place Jussieu, 75252 Paris Cedex 05, France} 

\address{$^\P$  Johann Radon Institute for Computational
 and Applied Mathematics, RICAM, Altenberger Strasse 69, A-4040 Linz,  Austria}

\vskip .2cm 

\begin{abstract}
  We show that the unresolved examples of Christol's conjecture
  $ \, _3F_{2}\left([2/9,5/9,8/9],[2/3,1],x\right)$ and
  $_3F_{2}\left([1/9,4/9,7/9],[1/3,1],x\right)$, are indeed diagonals of rational
  functions. We also show that other $\, _3F_2$ and $\, _4F_3$ unresolved examples of Christol's conjecture 
  are  diagonals of rational functions. Finally we give
  two arguments that show that it is likely that the  $\, _3F_2([1/9, 4/9, 5/9], \, [1/3,1], \, 27 \cdot x)$ function is a diagonal
of a rational function.
\end{abstract}

\vskip .1cm


\vskip .4cm


\noindent {\bf AMS Classification scheme numbers}: 33C05, 33C20, 33F10, 68W30

\vskip .2cm


\section{Introduction}
\label{Introduction}

There is a plethora of multiple integrals in physics:
Feynman integrals, lattice Green functions, the summands of the
magnetic susceptibility of the 2D Ising
model~\cite{Maillard2009,Maillard2009bis}, that have very
specific mathematical properties. These functions
are D-finite, i.e., solutions of linear differential
operators with polynomial coefficients,
and have series expansions with {\em integer coefficients}. Furthermore,
it was also shown that the linear differential operators
annihilating the summands of the magnetic susceptibility of the
Ising model $ \, \tilde{\chi}(n)$, verify the specific property of being
Fuchsian\footnote[1]{Denoting by $\, \theta$ the homogeneous
  derivative $\, x \cdot \, {{d} \over {dx}}$,
  the degrees of all the polynomial terms of the Fuchsian
  linear differential operator
  $\, \sum_n \, P_n(x) \cdot \, \theta^n$ are equal.} operators:
the critical
exponents of all their singularities are given by rational numbers,
and their Wronskians are $N$-th roots of
rational functions~\cite{Maillard2009,Maillard2009bis}. It was
also shown that the  $\tilde{\chi}(n)$ functions are solutions of
globally nilpotent operators~\cite{Maillard2009b}, and that
they ``come from geometry'' being G-operators~\cite{Andre1989}.

The unifying scheme behind these seemingly sparse properties
is that these functions are {\em diagonals of
  rational functions}~\cite{Short,Big}. It was shown for example
in~\cite{Big}, that if summands of the magnetic susceptibility
$ \, \tilde{\chi}(n)$ for any $n$ have an integer coefficient
series expansion reducing to algebraic series
modulo any prime, it is because they are diagonals of
rational functions for any integer~$n$. In fact many problems in
mathematical physics involving $n$-fold integrals, could be
interpreted in terms of diagonals of  algebraic or
rational\footnote[5]{Any diagonal of an algebraic
function in $n$ variables can be rewritten as the diagonal
of a rational function in $\, 2n$ variables: see~\cite{DenefLipshitz87}.}
functions\footnote[1]{See~\cite{Short,Big} p.26 and p.58.}.

Gilles Christol has shown in~\cite{christoldemo} that
for every rational function, its diagonal~$f(x)$ has the following properties:

\begin{itemize}
\item it is globally bounded: there exist integers $c$ and $d$
  in $\mathbb{N}^{*}$, such that $d \, f(c \, x) \in \mathbb{Z}[[x]] \, $
  and $f(x)$ has a radius of convergence that is non-zero in $\mathbb{C}$.
\item D-finite: there exists a differential operator
  $ \, L \in \mathbb{Z}[x][\frac{d}{dx}]$, with $L \, \neq 0$,
  such that $L(f) \, = 0$. 
\end{itemize}
Christol conjectured in~\cite{hypChristol} that every series
that verifies these two properties {\em is the diagonal of a rational function}.
In this paper~\cite{hypChristol} Christol came up with an unresolved example
to his conjecture~\cite{hypChristol}, and a longer list was generated
by Christol and his co-authors in 2012 in~\cite{Big}. In this paper
we show that two of the unresolved examples of the conjecture given in~\cite{Big}
on page 58, namely the
$ \, _3F_{2}\left([2/9,5/9,8/9], \, [2/3,1], \, 3^6 \cdot x\right)$
and $ \, _3F_{2}\left([1/9,4/9,7/9], \, [1/3,1], \, 3^6 \cdot x\right)$
are indeed diagonals of rational functions and provide a generalization of this result.

\vskip .2cm 
\vskip .2cm 

\section{Recalls on diagonals of rational functions and on Christol's conjecture}
\label{Recalls}

\subsection{Definition of the diagonal of a rational function}
The diagonal of a rational function in $\, n$ variables
$\,{\cal R}(x_1, \ldots, x_n)\, = \,\, $
${\cal P}(x_1, \ldots, x_n)/{\cal Q}(x_1, \ldots, x_n)$, 
where $ {\cal P}, {\cal Q} \, \in \mathbb{Q}[x_1, \ldots, x_n]$
such that $\, {\cal Q}(0, \ldots, 0) \neq 0$, is defined through its
 multi-Taylor expansion around $(0, \, \ldots, 0)$: 
\begin{eqnarray}
\label{defdiag}
\hspace{-0.90in}&&\quad \quad   \, \, \,
{\cal R}\Bigl(x_1, \, \ldots, \, x_n \Bigr)
\, \, \,\, = \, \,\, \,\sum_{m_1 \, = \, 0}^{\infty}
 \, \cdots \, \sum_{m_n\, = \, 0}^{\infty} 
 \,R_{m_1,  \, \ldots, \, m_n}
\cdot  \, x_1^{m_1} \,\,  \ldots \,\, x_n^{m_n}, 
\end{eqnarray}
as the series in  {\em one variable} $\, x$:
\begin{eqnarray}
\label{defdiag2}
\hspace{-0.7in}&&\quad \quad   \quad 
Diag\Bigl({\cal R}\Bigl(x_1, \, \ldots, \, x_n \Bigr)\Bigr)
\, \, \,\, = \, \,  \, \quad  \sum_{m \, = \, 0}^{\infty}
 \,R_{m, \, m, \, \ldots, \, m} \cdot \, x^{m}.
\end{eqnarray}

\vskip .1cm

\subsection{Hadamard product of algebraic functions and Christol's conjecture}

Recall that the Hadamard product of two series
$ \, f(x) = \, \sum_{n=0}^{\infty} \alpha_n \cdot \, x^n \, $ and
$ \, g(x) = \,  \sum_{n=0}^{\infty} \beta_n \cdot \,  x^n \, $ is given by:
 \begin{eqnarray}
\hspace{-0.90in}&& \quad \quad \quad  \quad \quad  \quad \quad \quad \, \, \,
f(x) \,  \star \,  g(x)  \,\,  \,
=  \, \,\,  \,  \sum_{n=0}^{\infty} \, \alpha_n \cdot \beta_n \cdot \,  x^n.
\end{eqnarray}

Hypergeometric series of the form
$ \, _pF_{p-1}([a_1, \, \ldots, \, a_p], \, [ b_1 \, \ldots, \, b_{p-1}], \, \, x) $
of height
$ \, h \, =  \,h(a_1, \, \ldots, \, a_p, \, b_1 \, \ldots, \, b_{p-1})$,
where the height $ \, h$ is given by:
\begin{eqnarray}
  \hspace{-0.90in}&&\quad \quad \quad \quad \, \, \,
h  =  \,  \, \,\,
\# \{ 1 \leq j \leq p \, \mid \, b_j \in \mathbb{Z} \}
\, \, - \# \{ 1 \leq j \leq p \, \mid \, a_j \in \mathbb{Z} \}
 \end{eqnarray}
 with $ \, b_p\, = \,  1$, can be written\footnote[2]{See~\cite{hypChristol} p.15.}
as the Hadamard product of $ \,h$ globally bounded\footnote[8]{Globally bounded
series can be recast into series with {\em integer} coefficients~\cite{Short,Big}.} series
of height $ \, 1$, were shown to verify Christol's conjecture. For example, the globally bounded
hypergeometric series $ \, _3F_2([1/3, \,1/3, \,1/3], \,[1, \,1],x)$ has height $3$, and it can be written as the Hadamard product of three algebraic functions\footnote[1]{Diagonals are closed under the Hadamard product: if
  two series are diagonals of rational functions, their
  Hadamard product is also a diagonal of a rational function.}:
\begin{eqnarray}
  \label{Diag1over3}
 \hspace{-0.98in}&& \,  \,  \, \,  \, \, \,
 _3F_2([{{1} \over {3}}, \, {{1} \over {3}}, \, {{1} \over {3}}], \, [1, \,1], \, x)
 \,\, = \, \, \, \,   (1-x)^{-1/3} \star (1-x)^{-1/3} \star (1-x)^{-1/3},
\end{eqnarray}
and can thus be written as the diagonal of the algebraic function in three variables:
\begin{eqnarray}
  (1-x)^{-1/3} \cdot (1-y)^{-1/3} \cdot (1-z)^{-1/3}.
\end{eqnarray}
Unlike the case of $ \, _3F_2([1/3, \,1/3, \,1/3], \,[1, \,1],x) \, $,
the hypergeometric functions
$ \, _3F_{2}([2/9,5/9,8/9], \, [2/3,1], \,   x)$ and
$ \, _3F_2([1/9,4/9,7/9], \, [1/3,1], \,  x)$,
while being {\em globally bounded} functions~\cite{Christol2},
were constructed in a way that avoids them
being written as ``simple'' Hadamard products of algebraic functions.

Note that a $_pF_{p-1}$ hypergeometric function is globally bounded without restrictions if all its parameters in the denominator are integers, while a $_pF_{p-1}$ hypergeometric function can be shown to be globally bounded in general, by looking at Landau functions as explained in the work of Christol~\cite{hypChristol}. Furthermore, Beukers and Heckman have shown in~\cite{Beukers}, that $_pF_{p-1}$ hypergeometric functions of height zero that are globally bounded are algebraic functions.

\vskip .1cm

\subsection{Unresolved examples to the conjecture}
\label{candidate}

Generalized hypergeometric functions with regular singularities $ \, _pF_{p-1}$
are a simple and natural testing ground for Christol's conjecture.

All $ \, \, _2F_1([a,b],[c],\,x)\,$  hypergeometric series
that are globally bounded are diagonals of rational functions.
There are two cases that fall into this category:

\begin{itemize}
\item Either the parameter $c \,$ is an integer,
and the $ \, _2F_1$ function can be written as the
Hadamard product of two $ \, _1F_0$ functions, which are algebraic functions,
and thus are diagonals of rational functions by Furstenberg's~\cite{Furstenberg}
theorem\footnote[5]{Furstenberg's theorem states that any algebraic function is the
diagonal of a rational function in two variables.}.  
\item The parameter $ \, c$ is not an integer, and in this case the $ \, _2F_1$ function
is a diagonal of a rational function
if and only if it is an algebraic function\footnote[9]{The only
  $ \, _2F_1$ hypergeometric functions that are globally bounded
  with $ \, c \in  \, \mathbb{Q}$ are the algebraic ones:
they are the ones appearing in the list of Schwarz~\cite{Schwarz}.}
(and consequently  its series is globally bounded).
\end{itemize}
Moving on to $\,  _3F_2$ hypergeometric functions, one can ask when is a $ \, _3F_2$
hypergeometric function a diagonal of a rational function?
\begin{itemize}
\item When the parameters $ \, d$ and $ \, e$ in $ \, _3F_2([a,b,c], \, [d,e], \, x)$
are integers,  because, in this case, the  $ \, _3F_2$ can be written as the
Hadamard product of three $ \, _1F_0$ algebraic functions,
and is thus the diagonal of a rational function, by the closure
  of diagonals under the Hadamard product.
\item When the parameters $d$ and $e$ in $ \, _3F_2([a,b,c], \, [d,e], \, x)$
 are rational numbers but not integers,
 because in this case the $ \, _3F_2$ is algebraic, and is thus
 a diagonal by Furstenberg's theorem.
\end{itemize}
Hence the interesting case occurs when {\em only one} of the two
parameters $ \, d$ or $ \, e$ is rational,
and the other is an integer. But even in this case, a lot of
the $\, _3F_2$ functions are easily seen to be diagonals of a rational function.
Suppose that a $ \, _3F_2([a,b,c], \, [1,e], \, x)$ is
globally bounded, with the parameters
$ \, a, \, b, \, c, \, e \in \,  \mathbb{Q} \setminus \mathbb{Z}$,
then there are six ways to write the $  \, _3F_2([a,b,c],[1,e],x)$
function as the diagonal of a rational function. This corresponds to the
six ways to write the $  \, _3F_2([a,b,c],[1,e],x)$ as a Hadamard product
of hypergeometric functions:
\begin{itemize}
\item  $_2F_1([a,b],\, [e],\, x) \star \, _1F_0([c], \, x)$
\item $_2F_1([a,c],\, [e],\, x) \star \, _1F_0([b], \, x)$
\item $_2F_1([b,c],\, [e],\, x) \star \, _1F_0([a], \, x)$
\item  $_2F_1([a,b],\, [1],\, x) \star  \, _2F_1([c,1],\, [e],\, x)$
\item   $_2F_1([a,c],\, [1],\, x) \star  \, _2F_1([b,1],\, [e],\, x)$
\item   $_2F_1([b,c],\, [1],\, x) \star  \, _2F_1([a,1],\, [e],\, x)$ 
\end{itemize}
Now $ \, _1F_0([c],\, x) \, $ and $_2F_1([a,b],[1],x) \, $
are diagonals of  rational functions
by what we have said above. Then $ \, _3F_2([a,b,c], \,[1,e],\, x)$
is a diagonal of  rational functions  if  $\, _2F_1([c,1], \,[e],\, x )\, $
or $  \, _2F_1([a,b],[e],\, x)$\footnote[9]{Instead of
$\, _2F_1([c,1],[e],\, x )$, or one could take any one of the three
  permuted versions: $\, _2F_1([b,1], \, [e],\, x )\, $, etc.},
are diagonals of rational functions, i.e.
if and only if they are algebraic functions,
since $ \, e  \, \in \, \mathbb{Q} \setminus \mathbb{Z}$. Now
$_2F_1([c,1],[e],x)$ cannot be an algebraic functions for $ \, e \in  \, \mathbb{Q}$
by Goursat~\cite{Goursat1881}. Hence if one of
$ \, _2F_1([a,b], \, [e], \, x)$,
$ \, _2F_1([b,c], \, [e], \, x)$, or $ \, _2F_1([a,c], \, [e], \, x)$
is an algebraic function, then  $ \, _3F_2([a,b,c],[1,e], \, x)$
is the diagonal of a rational function.
Now taking the two examples given in~\cite{Big} or by G. Christol
in~\cite{Christol2014}\footnote[8]{Appendix F p.58 of~\cite{Big},
  also see~\cite{Christol2014} p.19.}
that we are looking at here, we see that neither
$  \, _3F_{2}\left([2/9,5/9,8/9],[2/3,1], \, x \right)$,
nor $\,  _3F_{2}\left([1/9,4/9,7/9],[2/3,1], \,   x \right)$,
can be obtained as diagonals of rational functions through Hadamard
products\footnote[1]{See~\cite{rock} for a proof that $\,  _3F_{2}\left([1/9,4/9,7/9],[2/3,1], \,   x \right)$
  cannot be written as a Hadamard product.} since the three
$\, _2F_1$ hypergeometric series are {\em not}
globally bounded\footnote[5]{One can see this experimentally by taking
  the series expansion of any of the Gauss hypergeometric functions:
  the prime numbers in the denominators of the coefficients grow continuously.}:
\begin{eqnarray}
\label{3F2Christola31}
  \hspace{-0.98in}&& \quad \, \, 
 _2F_1\Bigl([{{2} \over {9}},\, {{5} \over {9}}], \, [{{2} \over {3}}], \, \, \,\,  x\Bigr),
\quad \quad
_2F_1\Bigl([{{2} \over {9}}, \, {{8} \over {9}}], \, [{{2} \over {3}}], \, \, \,\,  x\Bigr),
\quad \quad
_2F_1\Bigl([{{5} \over {9}}, \, {{8} \over {9}}], \, [{{2} \over {3}}], \, \, \,\,  x\Bigr),
\end{eqnarray}
and nor are the $\, _2F_1$ hypergeometric series:
\begin{eqnarray}
\label{3F2Christolb31}
  \hspace{-0.98in}&& \quad \, \, 
_2F_1\Bigl([{{1} \over {9}},\, {{4} \over {9}}], \, [{{1} \over {3}}], \, \, \,\,  x\Bigr),
\quad \quad
_2F_1\Bigl([{{4} \over {9}}, \, {{7} \over {9}}], \, [{{1} \over {3}}], \, \, \,\,  x\Bigr),
\quad \quad
_2F_1\Bigl([{{1} \over {9}}, \, {{7} \over {9}}], \, [{{1} \over {3}}], \, \, \,\,  x\Bigr).
\end{eqnarray}

\vskip .2cm 

\section{The main results}
\label{results}

The globally bounded $ \, _3F_{2}$ hypergeometric series
\begin{eqnarray}
  \label{3F2two}
 \hspace{-0.98in}&& \quad \quad \quad 
 _3F_2\Bigl([{{2} \over {9}}, \, {{5} \over {9}}, \, {{8} \over {9}}],
 \, [{{2} \over {3}}, \,1],  \,  \, 27 \cdot \, x \Bigr),
\quad \quad
_3F_2\Bigl([{{1} \over {9}}, \, {{4} \over {9}}, \, {{7} \over {9}}],
\, [{{1} \over {3}}, \, 1], \, 27 \cdot x \Bigr) 
\end{eqnarray}
are\footnote[2]{The operators annihilating the two hypergeometric
  functions (\ref{3F2two}) are adjoint of each other.}
respectively {\em the diagonals of the two algebraic functions}
\begin{eqnarray}
\label{11}
\hspace{-0.98in}&& \quad \quad \quad  \quad 
 _3F_{2}\Bigl([{{2} \over {9}},  \, {{5} \over {9}}, \, {{8} \over {9}}],
 \, [{{2} \over {3}}, \, 1],  \, \, 27 \cdot x\Bigr)
 \, \,\,\,  = \, \, \,\,\,   Diag\Bigl(\frac{(1-x-y)^{1/3}}{1-x-y-z}\Bigr),  
\end{eqnarray}
and
\begin{eqnarray}
  \label{12}
 \hspace{-0.98in}&& \quad  \quad \quad \quad 
 _3F_{2}\Bigl([{{1} \over {9}}, \, {{4} \over {9}}, \,  {{7} \over {9}}],
 \, [{{1} \over {3}}, \, 1], \,\,  27 \cdot x \Bigr)
 \, \, \,\,  = \, \, \,\,\,   Diag\Bigl( \frac{(1\, -x-y)^{2/3}}{1 \, -x-y-z}\Bigr). 
\end{eqnarray}
These two hypergeometric series\footnote[9]{The hypergeometric
  function  $\, _3F_2([2/9,5/9,8/9], \, [2/3,1],  27  x)$
  can be rewritten as the Hadamard product
  $_2F_1\Bigl([{{2} \over {9}}, \, {{5} \over {9}}], \, [{{2} \over {3}}],
  \, \, \,\,  27 \,x\Bigr)\,  \,  \star \,  \,  (1 \, -x)^{-8/9} \, $
  with  $\, _2F_1\Bigl([{{2} \over {9}}, \, {{5} \over {9}}],
  \, [{{2} \over {3}}], \, 27 \,x\Bigr)$ being 
  associated with a {\em Shimura curve}~\cite{Heun}. For more details
  please refer to \ref{AppendixA}.}
(\ref{3F2two}) can be recast into series with {\em integer} coefficients 
\begin{eqnarray}
\hspace{-0.98in}&& \, \, \, \, \, 
_3F_2\Bigl([{{2} \over {9}},  \, {{5} \over {9}}, \, {{8} \over {9}}],
                   \, [{{2} \over {3}}, \, 1], \, \,  3^6 \cdot\,  x\Bigr)
\, \, = \,\,\,  1 \,+  120 x \,+ 47124  x^2 \,+ 23483460 x^3 \, + \,\cdots \, ,
\end{eqnarray}
and
\begin{eqnarray}
 \hspace{-0.98in}&& \, \, \, \, \, 
 _3F_2\Bigl([{{1} \over {9}}, \, {{4} \over {9}}, \,  {{7} \over {9}}],
 \, [{{1} \over {3}}, \,1], \, \, 3^6 \cdot x)
\, \,= \,\,\, 1 \,\, +84 x \,+ 32760 x^2\,+ 16302000 x^3\,\,\, + \,\cdots \, 
\end{eqnarray}
Now Denef and Lipshitz in \cite{DenefLipshitz87} show that any power series in
$ \, \mathbb{Q} \left[ \left[x_1, \ldots, x_n  \right] \right]$,
algebraic over $ \, \mathbb{Q} (x_1, \ldots, x_n )$,
is the diagonal of a rational function in $2n$ variables, and
they give an algorithm to build this rational function. This means
that we can construct the rational functions, whose corresponding diagonals
are the $ \, _3F_{2}([2/9,5/9,8/9],[2/3,1], \, 27 \cdot x)$ and the
$ \, _3F_{2}([1/9,4/9,7/9],[1/3,1], \, 27 \cdot x)$
functions. We recall the algorithm of Denef and Lipshitz
and apply it to the algebraic function $(1-x-y)^{1/3}/(1-x-y-z)$
in the first subsection below, and then we give the rational function
and a generalization of the result in the second subsection. Finally,
we give a second proof of the general result using binomial sums. 

\vskip .2cm

\subsection{From diagonals of algebraic functions to diagonals of rational  functions: Denef and Lipshitz}

We explain a method which, for a given algebraic power series in $ \,n$ variables,
constructs a rational function in $\, 2n$ variables whose diagonal equals
the diagonal of the given algebraic series. Moreover, the partial diagonal of that $2n$-variable rational function, with respect to the pairs of variables $(x_1,x_{n+1}), \dots, (x_{n-1},x_{2n})$, yields the original $n$-variable algebraic power series. The method is described in the paper by Denef and Lipshitz~\cite{DenefLipshitz87} in the proof of their Theorem~6.2. As a running example we use the three-variable algebraic function
\begin{eqnarray}
  \hspace{-0.98in}&&  \quad  \quad  \quad  \quad  \quad  \quad  \quad   \quad  \quad 
 f(x, \, y, \,z)  \, \,  =  \,  \, \, \frac{(1-x-y)^{1/3}}{1-x-y-z},
\end{eqnarray}
whose multi-Taylor series expansion at~$0$ is actually a power series in the three
variables $\, x,y,z$:
\begin{eqnarray}
  \hspace{-0.98in}&&  \,  \,  
  f(x, \,y, \,z) \,  \,  =  \, \, \,\, 
 1 \, \, \,  + \,  {{2} \over {3}} \, x  \, + \,  {{2} \over {3}} \, y  \, + z  \,
 + \, {{10} \over {9}} \, x y  \, + \, {{5} \over {3}} \,  x z \,
 + \, {{5} \over {3}} \, y z \,  + \, {{40} \over {9}}  \,  x y z \, \,  \,  + \dots
\end{eqnarray}
Note that the minimal polynomial of $f$ is given by
\begin{eqnarray}
  \hspace{-0.98in}&& \quad  \quad \quad \quad \quad \quad 
  p(x,y,z,f) \, \,\,   = \, \, \,
  \Bigl((x + y + z - 1) \cdot \, f\Bigr)^3 \,\,\,   + 1 \,\, - x - y.
\end{eqnarray}
Denef and Lipshitz's theorem is formulated for \'{e}tale extensions, which
basically means that the partial derivative (w.r.t.~$f$) of the minimal
polynomial has a nonzero constant coefficient at~$0$. Clearly, the above
polynomial $\, p(x,y,z,f)$ does not meet this criterion.  However,
by considering $ \, \tilde{f}= \, f\, -1$, i.e. by removing the constant
term of~$f$, we can achieve an \'{e}tale extension. The minimal
polynomial then reads
\begin{eqnarray}
 \hspace{-0.98in}&& \quad  \quad \quad \quad \quad \quad 
\tilde{p}(x,y,z,f)\,\, \,  = \, \, \,
\Bigl((x + y + z - 1) \cdot (f+1) \Bigr)^3 \, \,  \, + 1 \,\, - x - y.
\end{eqnarray}
Indeed, $\frac{\partial\tilde{p}}{\partial f}(0,0,0,0)=-3 \neq 0$.  According to
the proof of Theorem~6.2~(i) in~\cite{DenefLipshitz87}, the rational
function
\begin{eqnarray}
  \hspace{-0.98in}&& \quad \quad \quad \quad \quad \quad \quad 
 \tilde{r}(x,y,z,f) \, \,\,  = \,\,\,
 f^2 \cdot \frac{\frac{\partial \tilde{p}}{\partial f}(xf,\, yf,\,zf, \, f)}{
   \tilde{p}(xf, \,yf, \,zf,\, f)}
\end{eqnarray}
has the property that
$ \, {\cal D}\bigl(\tilde{r}(x,\,y,\, z,\, f)\bigr) \, =\,\,\tilde{f}(x,y,z)$,
and hence $ \,  {\cal D}\bigl(r(x,y,z,f)\bigr) = \, f(x,y,z)\,$
for $\, r(x, \,y, \,z,\,f) \, = \, \, \tilde{r}(x,\,y,\,z,f)\,\, +1$. Here
the operator $\, {\cal D}$ denotes
a special kind of ``diagonalization'' with respect to the last variable: for
\begin{eqnarray}
  \hspace{-0.98in}&& \quad  \quad  \quad  \quad \quad 
f(x_1,\, \dots,\, x_n,\, y) \, \,  \, = \,\, \, \, 
\sum  \, a_{i_1,\dots,i_n,j} \cdot \,  x_1^{i_1} \cdots  \, x_n^{i_n} \, y^j, 
\end{eqnarray}
one defines
\begin{eqnarray}
\hspace{-0.98in}&& \quad \quad \quad \quad \quad 
  {\cal D}\bigl(f(x_1, \, \dots, \, x_n,y)\bigr)
  \,\,\, = \,\,
\sum_{j=i_1+\dots+i_n} \,  a_{i_1,\dots,i_n,j} \cdot \, x_1^{i_1} \cdots \,  x_n^{i_n}.
\end{eqnarray}
In our running example we obtain:
\begin{eqnarray}
 \hspace{-0.98in}&&  \quad \,  \quad 
r(x,\,y,\,z,\,f) \,\, \,  = \,\, \,\, 
\frac{3 \, f^2  \cdot \,  (f+1)^2 \cdot \,  (xf+yf+zf \,\,  -1)^3}{
  (f+1)^3  \cdot \,  (xf+yf+zf\, -1)^3 \,\,  -xf \, -yf \,\, \,  +1}
\, \, \, + \,1 .
\end{eqnarray}
In the second step, which is explained in the proof of Theorem~6.2(ii)
of~\cite{DenefLipshitz87}, one has to transform the rational function~$r$
(that has $ \, n+1$ variables) into another rational function (having $ \, 2n$
variables) such that its "true" (partial) diagonal gives the $n$-variable algebraic series~$f$. It
consists of a sequence of $ \, n -1$ elementary steps,
each of which is adding one more variable. In our example, we have to do the
following
\begin{eqnarray}
 \hspace{-0.98in}&& \quad 
r_1(x, \, y, \, z, \, u_1, \, v_1) \,\,   =\, \,\,
\frac{u_1 \cdot \,r(x, \, y, \, z, \, u_1)
  \, - v_1 \cdot \, r(x, \, y, \, z, \, v_1)}{u_1 \, -v_1},
 \\
\hspace{-0.98in}&& \quad 
r_2(x, \, y, \, z, \, u_1, \, u_2, \, v_2) \, \,  = \, \, \,
\frac{u_2 \cdot \,  r_1(x, \, y, \, z, \, u_1, \, u_2)
 \,  - v_2 \cdot \,  r_1(x, \, y, \, z, \, u_1, \, v_2)}{u_2 \, -v_2} \, ,
\nonumber
\end{eqnarray}
and obtain with $ \, r_2$ the desired rational function in {\em six} variables.

\vskip .2cm 

\subsection{Generalization of the previous result}

By the algorithm of Denef and Lipshitz given in the previous section,
  it is possible to show that the algebraic function 
\begin{eqnarray}
 \label{3F2Christol_alg}
 \hspace{-0.98in}&& \quad \quad  \quad \quad \quad \quad \quad \quad \quad \quad
 {{ (1 \, \, -x \, -y)^{a/b}  } \over {  1 \, -x \, - \, y -\, z  }},    
\end{eqnarray}
corresponds to the following rational function in {\em six} variables,
by taking the diagonal with respect to $ \, (x,u)$, $\, (y,v)$ and $\, (z,w)$:
\begin{eqnarray}
  \label{3F2ChristolbisA}
  \hspace{-0.98in}&& \quad 
 {{ a \cdot \, u^3 v \cdot \, (1\, -u x -u y -u z) \cdot \, (1 + \,u)^{a-1}
     \cdot \,  (1 \, -u x -u y -uz)^{a-1}} \over {
     (1 \, +u)^a  \cdot \, (1 \, -u x -u y -u z)^a
  \, -(1 \, -u x -u y)^b \cdot \,  (u \, -v) \cdot \,  (v \, -w)}}
 \nonumber \\
    \hspace{-0.98in}&& \quad  \, 
  -\, {{a  \cdot \, v^4  \cdot \, (1 \, -v x -v y -vz) \cdot \,  ((1 +\,v)
   \cdot \,  (1 \, -v x -v y -v z))^{a-1}} \over {
    (1 \, +v)^a \cdot \,  (1 \, -v x -v y -v z)^a
  \, -(1 \, -v x-vy)^b  \cdot \, (u \, -v)  (v \, -w)}}
  \\
\hspace{-0.98in}&& \quad  \,  \,  \, 
 - \, {{ a \cdot \,  u^3 w  \cdot \, (1 \, -u x-u y-u z) ((1 +\,u)
  \cdot \, (1 \, -u x -u y -u z))^{a-1}} \over  {
   (1 +\,u)^a \cdot \,  (1\, -u x -u y -u z)^a
  \, -(1 \, -u x -u y)^b \cdot \, (u \, -w) \cdot \, (v \, -w)}}
\nonumber  \\
  \hspace{-0.98in}&& \quad  \,  \,  \,  \,  \, 
  -\, {{a w^4  \cdot \, (1 \, -w x -w   y -w z)
      \cdot \,  (1 +\, w )^{a-1}  \cdot \, (1 \, -w x -w y -w z)^{a-1}} \over {
      (1 +\, w)^a \cdot \,  (1 \, -w x -w y -w z)^a
      \, - \, (1 \, -w x -w y)^b \cdot \, (u-w) \cdot \, (v-w)}}
 \, \,  \,  +1.
 \nonumber 
 \end{eqnarray}
The diagonal of the rational function (\ref{3F2ChristolbisA})
is annihilated by the linear differential operator of order three:
 \begin{eqnarray}
  \label{Christoloperator}
\hspace{-0.98in}&& \quad 
b^3\, x^2\, (1 -27 \,x) \cdot D_x^3
\, +b^2 \, x \,((27\,a\, -135\,b) \cdot \,x \, -a+3\,b) \cdot D_x^2
 \\
\hspace{-0.98in}&& \quad  \quad 
-b \cdot \,( (9\,a^2\, -63\,a\,b\, +114\,b^2) \cdot \,x \,  +a\,b-b^2) \cdot D_x
\, +(a-3\,b)\cdot (a-2 \,b) \cdot (a-b), \nonumber
\end{eqnarray}
 and can be expressed as the $\, _3F_2$ hypergeometric function 
\begin{eqnarray}
  \label{ChristolbisAhypergeom}
  \hspace{-0.98in}&& \quad \quad \quad \quad \quad
 _3F_2\Bigl([{\frac { 3\,a \, -b}{ 3\, a}}, \,  \,{\frac {2\,a \, -b}{3\, a}},
     \, \, {\frac {a \, -b}{3\, a}}],
 \, \, [{\frac { a \, -b}{a}}  , \, 1 ], \, \, 27  \cdot \, x \Bigr).
\end{eqnarray}

In particular, the two hypergeometric functions
$\, _3F_2([2/9,5/9,8/9], \, [2/3,1], \, 27 \cdot\,  x)$ and
$_3F_2([1/9,4/9,7/9],  \, [1/3,1], \, 27 \cdot \,  x)$
appearing in (\ref{3F2two}), correspond respectively
to the parameters $\, (a, \, b) \, = \, \, (1, \, 3)$,
and $\,  (a, \, b) \, = \, \, (2, \, 3)$ in the
algebraic function (\ref{3F2Christol_alg}). Other values of the
parameters $\,  (a, \, b)$ are not necessarily unresolved examples of Christol's conjecture.

For example if we consider the parameter values $\, a \, = \, \, 1 \, \, $
and  $\, b \, = \, \, 7 \, \, $, we see that the diagonal
of (\ref{3F2ChristolbisA}) is given by
the {\em globally bounded}\footnote[1]{$_3F_2\Bigl([{{2} \over { 7}}, \,  \, {{13} \over {21}},
    \, \, {{20} \over { 21}}],  \, \, [{{6} \over { 7}} , \, 1 ],
  \, \, 27 \cdot \, 7^3  \cdot \, x \Bigr) \, \, $
  is a series with {\em integer} coefficients.}
series (\ref{ChristolbisAhypergeom})
\begin{eqnarray}
  \label{ChristolbisAhypergeom}
  \hspace{-0.98in}&& 
  _3F_2\Bigl([{{2} \over { 7}}, \,  \, {{13} \over {21}}, \, \, {{20} \over { 21}}],
 \, \, [{{6} \over { 7}} , \, 1 ], \, \, 27 \, x \Bigr) \, \, = \, \, \, \,
 1 \, \, +{\frac {260 }{49}} \,x \,
 +{\frac {188190}{2401}}\, x^2  \,  \, \, + \, \, \cdots                 
\end{eqnarray}
with the $\, _2F_1$ series
\begin{eqnarray}
\label{ChristolbisAhypergeomextract}
\hspace{-0.98in}&& \, 
_2F_1\Bigl([{{13} \over {21}}, \, \, {{20} \over { 21}}],
\, \, [{{6} \over { 7}}], \, \, 27 \, x \Bigr),   \quad \, 
_2F_1\Bigl([{{2} \over { 7}},  \, \, {{20} \over { 21}}],
\, \, [{{6} \over { 7}}], \, \, 27 \, x \Bigr),   \quad \, 
_2F_1\Bigl([{{2} \over { 7}}, \,  \, {{13} \over {21}}],
\, \, [{{6} \over { 7}}], \, \, 27 \, x \Bigr),
 \nonumber 
\end{eqnarray}
being series that are {\em not} globally bounded. Hence
the hypergeometric series~(\ref{ChristolbisAhypergeom}) cannot be
easily written as a Hadamard product, as explained in Section~\ref{candidate}.

In contrast, for $\, a \, = \, \, 3 \, \, $ and $\, b \, = \, \, 4 \, \, $
the diagonal of (\ref{3F2ChristolbisA}) which is given by the
{\em globally bounded}\footnote[2]{$_3F_2\Bigl([{{3} \over {4}},
\,  \, {{5} \over {12}}, \, \, {{1} \over {12}}],  \, \, [{{1} \over { 4}} , \, 1 ],
  \, \, 1728 \cdot \, x \Bigr)\, \, $
  is a series with {\em integer} coefficients.}
series (\ref{ChristolbisAhypergeom8})
\begin{eqnarray}
  \label{ChristolbisAhypergeom8}
  \hspace{-0.98in}&& \quad
  _3F_2\Bigl([{{3} \over {4}}, \,  \, {{5} \over {12}}, \, \, {{1} \over {12}}],
 \, \, [{{1} \over {4}} , \, 1 ], \, \, 27 \, x \Bigr) \, \, = \, \, \, \,
 1 \,\,  \, +{\frac {45 }{16}} \, x\,  \, \, +{\frac {41769 }{1024}} \, x^2
  \,  \,\,  + \, \, \cdots      
\end{eqnarray}
with the $\, _2F_1$ series
\begin{eqnarray}
  \label{ChristolbisAhypergeomextract8}
  \hspace{-0.98in}&&  \quad \quad \quad \quad \quad \quad \quad \quad
  _2F_1\Bigl([{{5} \over { 12}}, \,  \, {{1} \over {12}}],
  \, \, [{{1} \over { 4}}], \, \, 27 \, x \Bigr),  
\end{eqnarray}
being a {\em globally bounded series}, which means that
it can be written as a diagonal using one of the procedures
given in Section~\ref{candidate}. We note that
algebraic functions close to the algebraic functions appearing
in (\ref{11}) and (\ref{12}), also give
$\, _3F_2$ or $\, _4F_3$
hypergeometric functions as their diagonals that are unresolved examples to Christol's conjecture:
\begin{eqnarray}
 \hspace{-0.98in}&& \quad \quad \quad 
 Diag\Bigl(\frac{(1-x-2 \, y)^{2/3}}{1-x-y-z}\Bigr)
 \, \, \,\,\, = \, \, \,\,\,
 _3F_{2}\Bigl([{{1} \over {9}},  \, {{4} \over {9}}, \, {{7} \over {9}}],
 \, [{{2} \over {3}}, \, 1],  \, \, 27 \cdot x \Bigr),
\end{eqnarray}
\begin{eqnarray}
 \hspace{-0.98in}&& \quad \quad \quad 
 Diag\Bigl(\frac{(1 \, -x-2 \, y)^{1/3}}{1\, -x-y-z}\Bigr)
 \, \, \,\,\, = \, \, \,\,\,
 _3F_{2}\Bigl([{{2} \over {9}},  \, {{5} \over {9}}, \, {{8} \over {9}}],
 \, [{{5} \over {6}}, \, 1],  \, \, 27 \cdot x\Bigr),
\end{eqnarray}
\begin{eqnarray}
 \hspace{-0.98in}&& \quad \quad \quad 
 Diag\Bigl(\frac{(1\, -x)^{1/3}}{1 \, -x-y-z}\Bigr)
 \, \, \,\,\, = \, \, \,\,\, _4F_{3}\Bigl([{{2} \over {9}},  \, {{5} \over {9}},
   \, {{8} \over {9}}, {{1} \over {2}}],
 \, [\, {{1} \over {3}}, \, {{5} \over {6}}, \, 1],  \, \, 27 \cdot x \Bigr),
\end{eqnarray}
\begin{eqnarray}
 \hspace{-0.98in}&& \quad \quad \quad 
 Diag\Bigl(\frac{(1 \, -x - y)^{1/3}}{1\, -x-z}\Bigr)
 \, \, \,\,\, = \, \, \,\,\,
 _4F_{3}\Bigl([{{2} \over {9}},  \, {{5} \over {9}}, \, {{8} \over {9}},\,
   {{-1} \over {3}}], \,
    [\, {{1} \over {3}}, \, {{5} \over {6}}, \, 1],  \, \, 27 \cdot x\Bigr).
\end{eqnarray}

\vskip .2cm

\subsection{Proof}
\label{Proof}
A computer algebra proof of this result can easily be obtained using
the creative telescoping program~\cite{Koutschan}: one
computes the operator (\ref{Christoloperator})
using the program~\cite{Koutschan}, and
verifies that this operator does annihilate the diagonal
of (\ref{3F2Christol_alg})\footnote[9]{One also needs to note
  that initial conditions have to be compared.}. Another longer
way to do it which we provide below, is through binomial sums.

The denominator of the algebraic function $ \, (1 \, -x-y)^{a/b}/(1 \, -x-y-z)$
can be expanded as a geometric series:
\begin{eqnarray}
\label{sum1}
 \hspace{-0.98in}&& \, \,\quad \quad \quad  \, \, \, 
 (1 \, -x-y-z)^{-1}\, \,  = \, \,\, 
 \sum_{n=0}^{\infty}\sum_{m=0}^{\infty}  {n \choose m} \cdot \,  (x+y)^{m} \,  z^{n-m}
 \nonumber \\
 \hspace{-0.98in}&& \, \, \, \, \,
 \quad \quad\quad\quad \quad\quad \quad\quad\quad\quad
  \, = \, \,\,   \sum_{n=0}^{\infty}\sum_{m=0}^{\infty}\sum_{l=0}^{\infty} \,
     {n \choose m} {m \choose l} \cdot \,  x^{l}\,y^{m-l}\,z^{n-m},
\end{eqnarray}
while the numerator can be written as the sum:
\begin{eqnarray}
\label{sum2}
 \hspace{-0.98in}&& \, 
 (1 \, -(x+y))^{a/b}\,\,  = \,\,\,
 \sum_{k=0}^{\infty} \frac{(-a/b)_{k}}{k!}\cdot   (x+y)^{k}
 \, = \, \,  \,
 \sum_{k=0}^{\infty} \sum_{j=0}^{k} \frac{(-a/b)_{k}}{k!}
 \cdot \,  {k \choose j}  x^{j}  y^{k-j}.
\end{eqnarray}
Multiplying these two sums \eqref{sum1} and~\eqref{sum2}
and re-indexing, we obtain: 
\begin{eqnarray}
\label{prodsums}
 \hspace{-0.98in}&& \, \quad  \, \, 
 \sum_{s=0}^{\infty}\sum_{t=0}^{\infty}\sum_{u=0}^{\infty} \,
 x^{s}\, y^{t}\, z^{u} \cdot \,\sum_{j=0}^{s} \sum_{k=0}^{ \infty} \,
 \frac{(-a/b)_{k}}{k!} \cdot \,  {k \choose j} \,
      {s+t+u-k \choose  s+t-k} \,  {s+t-k \choose  s-j} .
\end{eqnarray}
Now taking the coefficients corresponding to the diagonal in (\ref{prodsums}),
i.e. such that $s =  \, t = \,  u = \,  n$, we get:
\begin{eqnarray}
\label{b4chu}
 \hspace{-0.98in}&& \, \quad \quad \quad \quad \quad  \, \, \, \, 
 \sum_{j=0}^{n}\sum_{k=0}^{\infty}  \frac{(-a/b)_{k}}{k!} \cdot \,
   {k \choose j} \,  {3n-k \choose  2n-k} \,  {2n-k \choose  n-j}
  \nonumber \\
\hspace{-0.98in}&& \, \, \, \, \,
\quad \quad\quad\quad\,   \quad\quad\quad\quad\quad
\, = \, \,\,  \sum_{k=0}^{2n}   \frac{(-a/b)_{k}}{k!} \cdot \, {3n-k \choose  2n-k}
\cdot \, \sum_{j=0}^{n} \,
    {k \choose j} {2n-k \choose  n-j}.
\end{eqnarray}
Now recalling the Chu-Vandermonde identity which says that
$ {2n \choose n}\, = \,\, \sum_{j=0}^{n}\, {k \choose j}\, {2n-k \choose  n-j}$,
we find that (\ref{b4chu}) can be written as
\begin{eqnarray}
\label{afterchu}
\hspace{-0.98in}&& \, \, \, \, \,
 \quad\quad\quad\quad\quad \quad\quad \quad 
 S(n)\,  = \, \, {2n \choose n}  \cdot \,
 \sum_{k=0}^{2n}  \frac{(-a/b)_{k}}{k!} \cdot \,  {3n-k \choose  2n-k},
\end{eqnarray}
and by using a computer algebra tool like {\em Mathematica} or {\em Maple}
to simplify this sum into a closed form, from which we can read off
the hypergeometric function representation of the diagonal. More precisely,
we used creative telescoping (Zeilberger's algorithm)
to prove that (\ref{afterchu}) satisfies the first-order recurrence:
\begin{eqnarray}
  \label{recu}
  \hspace{-0.98in}&& \, \, \, \, \, \, \, \,
 \quad \quad   \quad \, \, \,
  (a\,-3 \, b \, -3\,b \,n ) \cdot (a\, -\, 2 \,b - \, 3\,b\,n)
  \cdot ( a\,- b\, - 3 \, b \, n ) \cdot \, S(n) \nonumber \\
\hspace{-0.98in}&& \, \, \, \, \,
 \quad\quad\quad\quad\, \,  \quad\quad \quad 
  = \, \, \, b^2 \cdot (n + 1 )^2 \cdot (a-  b  -  b \, n) \cdot \,  S(n+1). 
\end{eqnarray}
Together with the initial condition $S(0)=1$, we obtain the closed form
\begin{equation}
  S(n)\, = \,\, \frac{3^{3n} \cdot \, \bigl((b-a)/3b\bigr)_{n} \cdot\,
    \bigl((2b-a)/3b\bigr)_{n} \cdot\, \bigl((3b-a)/3b\bigr)_{n}}{
    \bigl((b-a)/b\bigr)_{n} \cdot \, \bigl(n!\bigr)^2}.
\end{equation}

\section{Telescopers of algebraic functions versus diagonals of algebraic functions }
\label{telescoperversusdiagonals}

The diagonal of a rational function and a solution of a
telescoper\footnote[1]{By ``telescoper'' of a rational function $\, R(x,y,z)$
we denote the output of the creative telescoping program~\cite{Koutschan},
applied to the transformed rational function $\, R(x/y,y/z,z)/(yz)$, which is
a differential operator that annihilates the diagonal of~$R$.} of a rational function
are close, yet distinct notions. A telescoper
annihilates an $n$-fold integral of a rational function over {\em all}
integration cycles\footnote[9]{Diagonals correspond {\em only} to
  {\em evanescent integration cycles} over rational functions.}. For
example the $\, _3F_2([a,b,c],[d,1],x)$ can be written through the well-known integral representation as:
\begin{eqnarray}
  \hspace{-0.98in}&&  \,  \quad \quad \quad \quad \quad \quad \quad
(1 - y)^{-1 - b + d} \cdot y^{b} \cdot ( 1 - x \cdot y^2)^{-a} \cdot (1 - z)^{-c},
\end{eqnarray}
with $\, a,b,c,d \in \mathbb{Q}$. Hence if one takes the parameters $\, a,b,c,d$
to have the values $a=1/9, \, b=4/9, \, c=7/9, \, d=1/3$,
one immediately obtains that the telescoper of the algebraic function
\begin{eqnarray}
\label{telescoper_alg}
  \hspace{-0.98in}&&  \quad \quad \quad \quad \quad \quad \quad \quad
\frac{y^{4/9}}{ (1-y)^{10/9} \cdot  (1-x \, y^2)^{1/9} \cdot  (1-z)^{7/9}} \, ,
\end{eqnarray}
admits as a solution the hypergeometric function
$_3F_{2}([{{1} \over {9}}, \, {{4} \over {9}}, \,  {{7} \over {9}}],
\, [{{1} \over {3}}, \, 1], \,\,  x)$. Yet the diagonal of the algebraic function 
(\ref{telescoper_alg}) is equal to zero, and it is through the
algebraic function (\ref{12}) that we were able to obtain it as the diagonal
of an algebraic function. Other $\, _3F_2 \, $ unresolved examples to Christol's
conjecture like~\cite{hypChristol,Christol2014}
\begin{eqnarray}
  \label{historic}
 \hspace{-0.98in}&& \quad \quad \quad \quad  \quad \quad \quad \
 _3F_2\Bigl([{{1} \over {9}},  \, {{4} \over {9}}, \, {{5} \over {9}}],
   \, [{{1} \over {3}}, \, 1],  \, \, 27 \cdot \,  x\Bigr),
\end{eqnarray}
were not obtained here as diagonals of a rational function, yet they
are solution of the telescoper of an algebraic function  and
can thus be seen as a {\em period of an algebraic variety over
a non-evanescent cycle}\footnote[2]{To be totally rigorous, one has
to consider the two certificates of the telescoping equation
see if that the integral of the derivatives of these two certificates
over that cycle are actually zero.}, but not necessarily as a diagonal
of an algebraic function (i.e. a period over an evanescent cycle). We
give two arguments in favour of the fact that the $ \, _3F_2$
hypergeometric function (\ref{historic}) is most probably a diagonal
of an algebraic function.

\vskip .1cm
 
\subsection{Diagonal: algebraic {\em mod} p}
\label{pFpmoins1}

If one expects $ \, _3F_2$ hypergeometric functions
like (\ref{historic}) to be diagonals of an
algebraic function, one should find~\cite{Short,Big} that
the corresponding series expansion
reduces to an algebraic series modulo any prime number $\, p$,
or power of a prime number $\, p^r$. In order to verify this fact
on (\ref{historic}) we look at the series expansion of
\begin{eqnarray}
\label{hyp}
\hspace{-0.96in}&&   \quad   \, \,
_3F_2\Bigl([{{1} \over {9}}, \, {{4} \over {9}}, \, {{5} \over {9}}],
\, [{{1} \over {3}}, \, 1],  \, \, 27^2 \cdot \, x\Bigr)
\, \, = \, \, \, \, \,
1 \,\, \, +60\,x \, \, \,+20475\,{x}^{2}\,  \, \,+9373650\,{x}^{3}
 \\
\hspace{-0.96in}&&  \quad \quad \quad  \quad
\, +4881796920\,{x}^{4} \, \, +2734407111744\,{x}^{5}
\, \, +1605040007778900\,{x}^{6}  \, \, + \,  \,  \cdots
\nonumber 
\end{eqnarray}
 which becomes modulo $\, 2$:
\begin{eqnarray}
\label{mod2}
\hspace{-0.96in}&&
_3F_2\Bigl([{{1} \over {9}}, \, {{4} \over {9}}, \, {{5} \over {9}}],
\, [{{1} \over {3}}, \, 1], \, \, 27^2 \cdot \, x\Bigr)
\, \,\,  = \, \, \, \,\,  \,
1 \,\,\, +{x}^{2} \, \,+{x}^{128} \, \, +{x}^{130}
\nonumber \\
\hspace{-0.96in}&&  \quad \quad  \, +{x}^{8192} \,+{x}^{8194} \,
+{x}^{8320} \,+{x}^{8322} \nonumber \\
\hspace{-0.96in}&&  \quad \quad   \,+{x}^{524288} \,+{x}^{524290}
\,+{x}^{524416} \,+{x}^{524418} \,  +{x}^{532480} \,
+{x}^{532482} \,+{x}^{532608} \,+{x}^{532610}
\nonumber \\
\hspace{-0.96in}&&  \quad \quad \quad \quad \quad \quad  \quad \quad
\,\, \,\,  +O({x}^{600000})
\nonumber \\
\hspace{-0.96in}&&  \quad 
\, \, = \, \, \, \, \, ( 1\, +{x}^{2})  \cdot \, (1\, + {x}^{128})
\cdot \, (1\, + {x}^{8192}) \cdot \,  (1 \, +{x}^{524288})
\,\,\,\,\, \,  \, \,  +O \left( {x}^{600000} \right).               
\end{eqnarray}
Straightforward guessing gives the infinite product formula
\begin{eqnarray}
\label{mod2guess}
\hspace{-0.96in}&& 
  F(x) \, \, = \, \,   \, ( 1\, +{x}^{2})  \cdot \, (1\, + {x}^{2^7})
\cdot \,  (1\, + {x}^{2^{13}}) \cdot \,  (1 \, +{x}^{2^{19}})
 \,    \, \cdots  \, \,  (1 \, + x^{2^{6\, n \, +1}})    \, \,  \,   \cdots  \,  
\end{eqnarray}
which is solution of
\begin{eqnarray}
\label{mod2guess}
\hspace{-0.96in}&& \quad \quad   \quad \quad      \quad \quad      
F(x)\, \, = \, \, \,   ( 1\, +{x}^{2})  \cdot \, F(x^{64})
\quad \quad  \quad  \quad mod. \quad  2, 
\end{eqnarray}
i.e. $\, _3F_2\Bigl([{{1} \over {9}}, \, {{4} \over {9}},
  \, {{5} \over {9}}], \, [{{1} \over {3}}, \, 1], 27^2 \cdot \, x\Bigr)$
is an algebraic function modulo $\, 2$ satisfying: 
\begin{eqnarray}
\label{mod2Eq}
\hspace{-0.96in}&& 
 \quad \quad     \quad \quad  \quad \quad
 F(x)\,\,  = \, \, \,   ( 1\, +{x}^{2})  \cdot \, F(x)^{64}
 \quad \quad  \quad \quad  mod. \quad  2
\end{eqnarray}
or:
 \begin{eqnarray}
\label{mod2Eq2}
\hspace{-0.96in}&&        
\quad \quad \quad \quad    \quad \quad      
( 1\, +{x}^{2}) \cdot \,  F(x)^{63} \,\, \,   = \, \,\, \,    1
\quad \quad  \quad \quad \quad \, \,  mod.  \quad 2.
 \end{eqnarray}
 Modulo $\, 3$ we have the following expansion
\begin{eqnarray}
\label{mod3}
\hspace{-0.96in}&& \quad \, \,   \quad  
 {{_3F_2\Bigl([{{1} \over {9}}, \, {{4} \over {9}}, \, {{5} \over {9}}],
\, [{{1} \over {3}}, \, 1], \, \, 27^2 \cdot \, x\Bigr) \, -1} \over {3}}
 \, \, = \, \, \, \, \,  2 \cdot \, F(x)
 \quad \quad  \quad \, \, mod. \quad  3, 
\end{eqnarray}
where:
\begin{eqnarray}
\label{mod3}
\hspace{-0.96in}&&\quad   \quad \quad  
 F(x)    \, \, \, = \, \,\,  \,\,
 x \, \, \,+{x}^{3}\,\,\,  +{x}^{9}\,  \,+{x}^{27}\,
 +{x}^{81}\,\, +{x}^{243}\,\, +{x}^{729}\,\, +{x}^{2187}\,\, +{x}^{6561}
 \nonumber \\
  \hspace{-0.96in}&&  \quad\quad \quad \quad \quad \quad  \quad \quad \quad 
\, +{x}^{19683}\,\, \, +{x}^{59049} \,\, \,\, +O \left( {x}^{60000} \right)               
\end{eqnarray}
which is solution of
\begin{eqnarray}
\label{mod3guess}
\hspace{-0.96in}&& \quad \quad  \quad  \quad \quad \quad      \quad \quad      
x \,  \, + \, \,  F(x^3) \, \, = \, \, \,  \, F(x)
\quad \quad  \quad \quad  \quad  mod. \quad  3, 
\end{eqnarray}
i.e. $\, F(x)$ is an algebraic function modulo $\, 3$: 
\begin{eqnarray}
\label{mod2Eq}
\hspace{-0.96in}&& 
 \quad \quad     \quad \quad  \quad \quad \quad  \quad
 x \, \,  + \, \,  F(x)^3 \, \, = \, \, \,  \, F(x)
 \quad \quad  \quad \quad  \quad  mod. \quad  3. 
\end{eqnarray}
Unlike the situation on the
$ \, _3F_2\bigl([{{1} \over {9}}, \, {{4} \over {9}}, \, {{7} \over {9}}],
\, [{{1} \over {3}}, \, 1],  \, \, 27^2 \cdot \, x\bigr)$ hypergeometric series,
it is less obvious how to obtain the
$_3F_2\bigl([{{1} \over {9}}, \, {{4} \over {9}}, \, {{5} \over {9}}],
\, [{{1} \over {3}}, \, 1],  \, \, 27^2 \cdot \, x\bigr)$
as the diagonal of a rational function. It is however possible
to obtain the solution of
$ \, _3F_2\Bigl([{{1} \over {9}}, \, {{4} \over {9}}, \, {{5} \over {9}}],
\, [{{1} \over {3}}, \, 1],  \, \, 27^2 \cdot \, x\Bigr)$,
as a telescoper of an algebraic function, and this solution is an algebraic function modulo p.

\subsection{A relation between $\, _3F_2([1/9, 4/9, 5/9], \, [1/3,1], \, 27 \cdot x)$
     and a $ \, _4F_3$ diagonal of an algebraic function}

The diagonal of the product of algebraic functions
\begin{eqnarray}
  \hspace{-0.96in}&& \, \quad \quad  \quad \quad  \quad \quad  \quad \quad 
 \frac{(1 - x - y)^{2/3}}{(1 - x - y - z)} \cdot \, (1 - w)^{-5/9}, 
\end{eqnarray}
is given by the  $ \, _4F_3$ hypergeometric function $\, {\cal H}$
which is the Hadamard product of
$\, _3F_2([1/9, 4/9, 7/9], \, [1/3,1], \, 27 \cdot x)$ and $\, (1 \, -x)^{-5/9}$: 
\begin{eqnarray}
  \label{4F3chris}
\hspace{-0.96in}&& \, \quad \quad \quad  \quad  \quad \, 
       {\cal H} \, = \, \,
       _4F_3\Bigl([{{1} \over {9}}, \,  {{4} \over {9}}, \,  {{5} \over {9}},
  \, {{7} \over {9}}], \, [{{1} \over {3}}, \,  1, \,  1], \, \,  27 \cdot x \Bigr)
\nonumber \\
\hspace{-0.96in}&& \, \quad \quad \quad \quad \quad \quad   \quad \quad \, 
\, \, = \, \, \,
 (1 \, -x)^{-5/9} \,  \star  \,\,
_3F_2\Bigl([{{1} \over {9}}, \,  {{4} \over {9}}, \, {{7} \over {9}}],
\, [{{1} \over {3}}, \, 1], \,  \, 27 \cdot x\Bigr)
\nonumber \\
 \hspace{-0.96in}&& \, \quad \quad \quad \quad \quad  \quad   \quad  \quad   
 \, \, = \, \, \, \, \, 
 Diag\Bigl(\frac{(1 - x - y)^{2/3}}{
   (1 - x - y - z)} \cdot \, (1 - w)^{-5/9}\Bigr).
\end{eqnarray}
This $\, _4F_3$ hypergeometric series  (\ref{4F3chris})
{\em can also} be written as the Hadamard product:
\begin{eqnarray}
  \label{4F3chrisHad}
  \hspace{-0.96in}&& \, \quad \quad \quad  \quad \quad \quad 
 {\cal H} \, = \, \,  (1 \, -x)^{-7/9} \,  \star  \,\, 
  _3F_2([{{1} \over {9}}, \,  {{4} \over {9}}, \, {{5} \over {9}}],
  \, [{{1} \over {3}}, \, 1], \,  \, 27 \cdot  \, x).
\end{eqnarray}
So even though we did not find a rational (or algebraic)
function whose diagonal
is given by  (\ref{historic}),
knowing that  $ \, _3F_2\Bigl([{{1} \over {9}}, \,  {{4} \over {9}},
  \, {{7} \over {9}}],[{{1} \over {3}}, \, 1], \,  \, 27 \cdot x\Bigr)$
is the diagonal of a rational function, we see that the Hadamard product
of (\ref{historic}) with a simple algebraic function   $\, (1 \, -x)^{-7/9}$
is actually a diagonal of an algebraic (or rational) function.
This suggests but does not prove, that
$ \, _3F_2\Bigl([{{1} \over {9}}, \,  {{4} \over {9}},
  \, {{5} \over {9}}],[{{1} \over {3}}, \, 1], \,  \, 27 \cdot x\Bigr)$ could also
be a diagonal of a rational function.

\vskip .1cm 
\vskip .1cm 

\section{Conclusion}
\label{Conclusion}

Because of the crucial role played by diagonals of rational functions
in physics~\cite{Short,Big}, Christol's conjecture is an important open problem.
We have shown that the hypergeometric series
$ \, _3F_{2}\left([2/9,5/9,8/9], \, [2/3,1], \,  x\right)$ and
$_3F_{2}\left([1/9,4/9,7/9], \, [1/3,1], \,  x\right)$ appearing
in~\cite{Big,Christol2014} are diagonals of rational functions. We did so
by first finding two algebraic functions whose diagonals were given
by these two hypergeometric functions, and through an algorithm outlined
in the paper~\cite{DenefLipshitz87}, we were able to recover the rational
functions whose diagonals are given by these two $\, _3F_{2}$
hypergeometric functions.

We were also able to give a generalization of
this result, and obtain other unresolved examples of Christol's conjecture
as diagonals of rational functions. Furthermore, even though we were
not able to write the $\, _3F_2([1/9, 4/9, 5/9], \, [1/3,1], \, 27 \cdot x)$
given by Christol in~\cite{hypChristol}, as a diagonal
of a rational function,
we gave two arguments that suggested that it was likely
to be so. More generally, we believe after writing the
$ \, _3F_{2}\left([2/9,5/9,8/9], \, [2/3,1], \,  x\right)$ and
$_3F_{2}\left([1/9,4/9,7/9], \, [1/3,1], \,  x\right)$ as diagonal
of rational functions, that it is likely that the other
$\, _3F_2$  unresolved examples of Christol's conjecture
are diagonals of rational functions.

\vskip .3cm

\vskip .3cm

{\bf Acknowledgments.}  
J-M. M. and Y.A. would like to thank G. Christol for many enlightening 
discussions on diagonals of rational functions. 
We would like to thank  A. Bostan 
for useful discussions on creative telescoping. Y.A. would like to thank A. Bostan for many explanations
on Christol's conjecture
including the details of \ref{candidate} in an
enlightening email correspondence~\cite{Bostan}.
J-M. M. would like to thank the School of Mathematics
and Statistics of Melbourne
University where part of this paper was written for hospitality. 
J-M. M. would like to thank A.J. Guttmann for many discussions on D-finite series. We would like to thank A.J. Guttmann for proof-reading the paper. 
 Y.A. would like to thank J. Voight for providing enlightening explanations on Shimura curves and for providing reference~\cite{Takeuchi1977}.
C.K. was supported by the Austrian Science Fund (FWF): P29467-N32.
Y. A. was supported by the Austrian Science Fund (FWF): F5011-N15.
Y. A. would like to thank the RICAM for hosting him on several occasions. Y.A. would like to thank Elaine and Rob for their hospitality in Linz, and for the great discussions he had with them there. Y.A. would like to thank his parents and his family for all their support and love. We thank the Research Institute for Symbolic Computation for giving us access to the RISC software packages. 

 \vskip .5cm 
\appendix

\section{Counterexamples and links with Shimura curves}
\label{AppendixA}
 The  Gauss hypergeometric function appearing on the left in (\ref{3F2two})
 \begin{eqnarray}
   \hspace{-0.98in}&& \quad  \quad \quad  \quad \quad  \quad  \quad   \quad  \quad  
   \, _3F_2\Bigl([{{2} \over {9}},\, {{5} \over {9}}, \, {{8} \over {9}}],
   \, [{{2} \over {3}}, 1], \, \, \,\,  27 \,x \Bigr)
\end{eqnarray}
 can be seen as the Hadamard product of a Gauss hypergeometric function
 and an algebraic function given by:
\begin{eqnarray}
 \hspace{-0.98in}&& \quad  \quad\quad  \quad \quad  \quad  \quad   \quad  \quad  
 \, _2F_1\Bigl([{{2} \over {9}}, \, {{5} \over {9}}],
 \, [{{2} \over {3}}], \, \, \,\,  27 \,x\Bigr)
\,  \,  \star \,  \,  (1 \, -x)^{-8/9}.
\end{eqnarray}
Now the  Gauss hypergeometric function
$\,\, _2F_1\Bigl([{{2} \over {9}},\, {{5} \over {9}}],
\, [{{2} \over {3}}], \, \, \,\,  27 \,x\Bigr) \,\,$
happens to be a hypergeometric function corresponding to an automorphic form
associated with a {\em Shimura curve}~\cite{Heun,Voight,Voight2}.
One has the identity:
\begin{eqnarray}
\label{3F2Christola31b7}
 \hspace{-0.98in}&& \quad  \quad  \quad
 _2F_1\Bigl([{{2} \over {9}},\, {{5} \over {9}}], \, [{{2} \over {3}}],
 \, \, \,\,  27 \, x\Bigr) \, \, \, = \, \, \, \, \,
  (1 \, -27\, x)^{-1/9} \cdot \, (1 \, -36\, x \, +216\, x^2)^{-1/18}  
\nonumber \\
  \hspace{-0.98in}&& \quad \quad \quad \quad \quad \quad \quad
 \, \, \times \,   \,
_2F_1\Bigl([{{1} \over {36}},\, {{19} \over {36}}], \, [{{8} \over {9}}],
 \, \, \,\,  -\, 1728 \cdot \,{\frac {{x}^{3} \cdot \, (1\, - 27\,x) }{
     (1 \, -36 \, x \, +216 \, x^2)^2}}\Bigr).
\end{eqnarray}
The  Gauss hypergeometric function
$\, _2F_1\Bigl([{{1} \over {36}},\, {{19} \over {36}}],\, [{{8} \over {9}}], \, \,x\Bigr)$
can be also expressed as:
\begin{eqnarray}
\label{3F2Christola31bc7}
  \hspace{-0.98in}&& \quad \, \,\,
 _2F_1\Bigl([{{1} \over {36}},\, {{19} \over {36}}], \, [{{8} \over {9}}], \, \,x\Bigr)
 \, \,\, = \, \,\, \, (1\, -x)^{-1/36} \cdot \,
 _2F_1\Bigl([{{1} \over {36}},\, {{13} \over {36}}],
 \, [{{8} \over {9}}], \, \, -\, {{x} \over {1\, - x}}\Bigr).
\end{eqnarray}
Now the  Gauss hypergeometric function
$\, _2F_1([{{1} \over {36}},\, {{13} \over {36}}], \, [{{8} \over {9}}], \, \,x) \, $
which occurs in p.14 of~\cite{Hoeij2014}, corresponds to
a hypergeometric function related to a {\em  Shimura curve} since
it has exponent differences\footnote[1]{See~\cite{Tu2013} p.10
  for a definition of exponent difference.}
$ \,(1/9, \,1/2, \, 1/3 )$, and these exponent differences are listed
in the exhaustive list of hypergeometric functions that are associated
with Shimura curves appearing in Table~1 of~\cite{Takeuchi1977}.
Other $\, _3F_2$ functions that are unresolved examples to Christol's conjecture
that we found to be related to $\, _2F_1$
hypergeometric functions related to {\em Shimura curves} are given by:
\begin{eqnarray}
\label{147morecidentitybis}
 \hspace{-0.96in}&&   \, \quad 
 _3F_2\Bigl([{{1} \over {9}}, \, {{4} \over {9}}, \, {{7} \over {9}}],
 \, [ {{4} \over {3}}, \, 1],  \, \,  \,3^6\,  x\Bigr) \, \, = \, \, \,
 (1\, -x)^{-1/9} \, \,  \star \,\,  _2F_1\Bigl([{{4} \over {9}},
   \, {{7} \over {9}}], \, [ {{4} \over {3}}],  \, \,  \,3^6\,  x\Bigr), 
\end{eqnarray}
\begin{eqnarray}
\label{reduce1}
 \hspace{-0.96in}&&   \, \quad 
 _3F_2\Bigl([{{2} \over {9}}, \, {{5} \over {9}}, \, {{7} \over {9}}],
 \, [ {{2} \over {3}}, \, 1],  \, \,  \,3^6\,  x\Bigr) \, \, = \, \, \,
 (1\, -x)^{-7/9} \, \,  \star \,\,  _2F_1\Bigl([{{2} \over {9}},
   \, {{5} \over {9}}], \, [ {{2} \over {3}}],  \, \,  \,3^6\,  x\Bigr), 
\end{eqnarray}
\begin{eqnarray}
\label{reduce2}
 \hspace{-0.96in}&&   \, \quad 
 _3F_2\Bigl([{{4} \over {9}}, \, {{5} \over {9}}, \, {{8} \over {9}}],
 \, [ {{2} \over {3}}, \, 1],  \, \,  \,3^3\,  x\Bigr) \, \, = \, \, \,
 (1\, -x)^{-8/9} \, \,  \star \,\,  _2F_1\Bigl([{{4} \over {9}},
   \, {{5} \over {9}}], \, [ {{2} \over {3}}],  \, \,  \,3^3\,  x\Bigr), 
\end{eqnarray}
\begin{eqnarray}
\label{reduce3}
 \hspace{-0.96in}&&   \, \quad 
 _3F_2\Bigl([{{1} \over {7}}, \, {{2} \over {7}}, \, {{4} \over {7}}],
 \, [ {{1} \over {2}}, \, 1],  \, \,  \, 7^4 \,  x\Bigr) \, \, = \, \, \,
 (1\, -x)^{-4/7} \, \,  \star \,\,
 _2F_1\Bigl([{{1} \over {7}},  \, {{2} \over {7}}],
 \, [ {{1} \over {2}}],  \, \,  \,7^4 \,  x\Bigr), 
\end{eqnarray}
Besides two hypergeometric functions, the
$\, _3F_2([{{2} \over {9}},\, {{5} \over {9}}, \, {{8} \over {9}}],
\, [{{2} \over {3}}, 1], \, \, \,\,  27 \,x)$ and
the  $\, _3F_2$ hypergeometric
$\, _3F_2([{{1} \over {9}}, \, {{4} \over {9}}, \, {{7} \over {9}}],
\, [ {{4} \over {3}}, \, 1],  \, \,  \,27\,  x)$,
and the three  globally bounded $\, _3F_2$ hypergeometric series
(\ref{reduce1}), (\ref{reduce2}) and  (\ref{reduce3}), 
we were not able to write the other examples given in this section
as a Hadamard product involving a $\, _2F_1$
hypergeometric function associated to a {\em Shimura curve}. In any case,
since the class of potential counterexamples formulated by Christol
is {\em infinite}, while the list of Shimura in  Table 1 of~\cite{Takeuchi1977}
is {\em finite}, a list of  $\, _3F_2$ functions both
related to {\em Shimura curves} and to Christol's
conjecture is bound to be finite. 

\vskip .3cm  

\vskip .3cm  

\vskip .3cm  



\section*{References}

\begin{thebibliography}{10}


\bibitem{Bostan} A. Bostan ``Re: "Trivialement diagonale"''.
  Message to Youssef Abdelaziz. 16 March 2019.
  
 \bibitem{Maillard2009} A. Bostan, S. Boukraa, A.J. Guttmann, S. Hassani,
   I. Jensen, J-M. Maillard and Zenine
  {\em High order Fuchsian equations for the square lattice Ising model: $\tilde{\chi}^{(5)}$}
  \newblock  J. Phys. {\bf A 42} 275209, (2009),
  \newblock \url{http://arxiv.org/abs/0904.1601}

\bibitem{Short}
  A. Bostan, S. Boukraa,  G. Christol, S. Hassani,  J-M. Maillard,
{\em Ising n-fold integrals as diagonals of rational functions and 
integrality of series expansions},  (2013), 
\newblock  J. Phys. {\bf A 46}: Math. Theor.  185202 (44 pages),     
\newblock  \url{http://arxiv.org/abs/1211.6645v2} 

\bibitem{Big}
  A. Bostan, S. Boukraa, G. Christol, S. Hassani and {J-M} Maillard
  {\em Ising $\,n$-fold integrals as diagonals of rational functions
    and integrality of series expansions: integrality versus modularity} {P}reprint,  (2012)
 \newblock  \url{http://arxiv.org/abs/1211.6031}
  
\bibitem{Maillard2009b} A. Bostan, S. Boukraa, S. Hassani, J-M. Maillard, J-A. Weil and N. Zenine
\newblock  {\em  Globally nilpotent differential operators and the square Ising model}
\newblock  J. Phys. {\bf A 42} 125206, (2009)
\newblock \url{http://arxiv.org/abs/0812.4931}

\bibitem{Goursat1881} E. Goursat,
  {\em Le{c}cons sur les s\'eries hyp\'erg\'eom\'etriques et sur quelques
    fonctions qui s'y rattachent}, Hermann (Paris), (1938).

\bibitem{Beukers} F. Beukers, G. Heckman
  {\em Monodromy for the hypergeometric function nFn-1}
   Inventiones mathematicae 95.2, 325-354 (2009)
  
\bibitem{Tu2013}   F-T. Tu,
  {\em Algebraic Transformations of hypergeometric functions arising from theory of Shimura curves}
  \newblock RIMS K\^oky\^uroku Bessatsu, {\bf B 44}, pp.223-245 (2013).

\bibitem{christoldemo} G. Christol,
  {\em Diagonales de fractions rationnelles et \'{e}quations de Picard-Fuchs},
  G.E.A.U. 12\`eme ann\'ee, 1984/85, n\textsuperscript{o} 13, 12p.

\bibitem{hypChristol}
  G. Christol, {\em Fonctions hyperg\'eom\'etriques born\'ees},  (2013),
  Groupe de travail d'analyse ultram\'etrique, {\bf 14} (1986-1987),
  Expos\'e No. 8, 16 pp.

\bibitem{Christol2014} G. Christol,
  {\em Fonctions hyperg\'eom\'etriques et diagonales de fonctions rationnelles. La fonction 
  $\, _3F_2([1/9, 4/9, 5/9], \, [1/3,1],  x)$ est-elle une diagonale?}
  Conference ``Journ\'ees holonomes'',
  Institut Fourier, Grenoble, February 2014. Slides available at
\url{http://www-fourier.ujf-grenoble.fr/~jroques/journeesholonomes.html}

\bibitem{Christol2} G. Christol
  { \em Globally bounded solutions of differential equations}
  \newblock  Analytic number theory (Tokyo, 1988) (Lecture Notes in Math. vol 1434)
  (Berlin: Springer) pp. 45-64
\newblock \url{http://dx.doi.org/10.1007/BFb0097124}

\bibitem{Schwarz} H.A. Schwarz 
  \newblock {\em Ueber diejenigen F\"alle, in welchen die Gaussische
    hypergeometrische Reihe eine algebraische Function ihres vierten Elementes darstellt}
  \newblock {Journal f\"ur die reine und angewandte Mathematik
    (Crelle's Journal), Volume 1873 (75) Jan 1, 1873, 44 pages}

 \bibitem{Koutschan} HolonomicFunctions Package version 1.7.1 (09-Oct-2013)
written by Christoph Koutschan, Copyright 2007-2013, Research Institute for 
Symbolic Computation (RISC),
Johannes Kepler University, Linz, Austria

\bibitem{Furstenberg} H. Furstenberg
  {\em Algebraic functions over finite fields}
   \newblock J. Algebra {\bf 7}, 271-277 (1969)
   \newblock \url{http://dx.doi.org/10.1016/0021-8693(67)90061-0}

 \bibitem{Koutschan} HolonomicFunctions Package version 1.7.1 (09-Oct-2013)
written by Christoph Koutschan, Copyright 2007-2013, Research Institute for 
Symbolic Computation (RISC),
Johannes Kepler University, Linz, Austria

\bibitem{DenefLipshitz87} J. Denef, L. Lipshitz,
  {\em Algebraic power series and diagonals}, 
  \newblock J. Number Theory {\bf 26} 46-67
  \newblock \url{http://dx.doi.org/10.1016/0022-314X(87)90095-3}

\bibitem{Voight}  J. Voight,
  {\em Shimura curves of genus at most two}, Math. Comp. 78 pp 1155-1172, (2009).

\bibitem{Voight2}  J. Voight,
  {\em Three lectures on Shimura curves}, 16th april (2006).

\bibitem{Takeuchi1977}  K. Takeuchi,
  {\em Commensurability classes of arithmetic triangle groups}, 1977

\bibitem{Hoeij2014}  M. van Hoeij and R. Vidu\~{n}as,
  \newblock {\em Belyi functions for hyperbolic Hypergeometric-to-Heun transformations},
  \newblock Journal of Algebra, Volume {\bf 441}, 1 November 2015, Pages 609-659, (2015)
  \url{arXiv:1212.3803v3} 

 \bibitem{Maillard2009bis} S. Boukraa, S. Hassani,
   I. Jensen, J-M. Maillard and N. Zenine 
  {\em High order Fuchsian equations for the square lattice Ising model: $\tilde{\chi}^{(6)}$}
  \newblock  J. Phys. {\bf A 43} 115201, (2010),
  \newblock \url{http://arxiv.org/abs/0912.4968}

  
\bibitem{rock} T. Rivoal, J. Roques {\em Hadamard products of algebraic functions},
  Journal of Number Theory, 145:579--603, (2014). 
  
\bibitem{Heun} Y. Abdelaziz, S. Boukraa,  C. Koutschan,  J-M. Maillard,
   \newblock {\em Heun  functions and diagonals of rational functions (unabridged version)},
    \url{arXiv:1910.10761v1}
   
\bibitem{Andre1989} Y. Andr\'e
  {\em  G-functions and geometry Aspects of Mathematics}
  \newblock E13 (Braunschweig: Friedr.Vieweg and Sohn) ISBN 3-528-06317-3

\vskip .3cm

\end{thebibliography}
\end{document}